\documentclass[11pt]{article}
\usepackage{amssymb}
\usepackage{amsmath}
\usepackage{graphicx}
\usepackage[utf8]{inputenc}
\usepackage{amsfonts}
\usepackage{cite}
\usepackage{graphicx}
\usepackage{xcolor}
\usepackage[utf8]{inputenc}
\usepackage{float}
\usepackage{hyperref}
\usepackage{tabularx}
\usepackage{breqn}
\usepackage[english]{babel}
\def\bc{\begin{center}}
	\def\ec{\end{center}}

\def\s2c{\vskip 2cm}
\def\bt{\begin{Theorem}}
	\def\et{\end{Theorem}}
\def\bd{\begin{Definition}}
	\def\ed{\end{Definition}}
\def\bl{\begin{Lemma}}
	\def\el{\end{Lemma}}
\def\be{\begin{Example}}
	\def\ee{\end{Example}}
\def\bcor{\begin{Corollary}}
	\def\ecor{\end{Corollary}}
\def\br{\begin{Remark}}
	\def\er{\end{Remark}}

\def\mysection{\setcounter{equation}{0}\section}
\newtheorem{Lemma}{Lemma}[section]
\newtheorem{Theorem}[Lemma]{Theorem}
\newtheorem{Definition}[Lemma]{Definition}
\newtheorem{Proposition}[Lemma]{Propostion}
\newtheorem{Corollary}[Lemma]{Corollary}
\newtheorem{Remark}[Lemma]{Remark}
\date{}
\title{Analysis of steady state solutions to an age structured SEQIR model with optimal vaccination }
\author {Manoj Kumar, Syed Abbas \\
School of Basic Sciences,\\
Indian Institute of Technology Mandi,\\	Kamand (H.P.) - 175005, India
	\\Email :  sabbas.iitk@gmail.com; abbas@iitmandi.ac.in}
	
\begin{document}
	\maketitle
	\author
	\noindent {\bf Abstract} :
  Quarantine of those individuals who are suspected of being infected is one of the intervention measures to contain the spread of an infectious disease. We propose an age structured SEQIR (S-Susceptible, E-Exposed, Q-Quarantine, I-Infected, R-Recovered) model with vaccination of susceptible and exposed individuals. Firstly without vaccination, basic reproduction number is derived by using the appearance of endemic steady state. With an appropriate Liapunov function, stability of disease free equilibrium point is checked. Vaccinating a population with reduction in number of infected individuals at minimal cost is considered as an optimization problem. We show that the vaccination strategy is concentrated on atmost three age classes.

\mysection{Introduction}
Introduction of age structure in classical SIR models allows us to model more realistic scenarios related to the disease progression, but it also produces complex dynamics and possibly a change in the behavior of the solutions. As it is also observed in the case of COVID-19, age plays a crucial role for the spread of disease. Adding age structure to simple SIR models also increases the complexity of the models and the qualitative analysis becomes more complicated as we need more advanced techniques to handle these models in a sophisticated manner. Kermack and A.G. McKendrick in a sequence of famous research papers \cite{RN57, RN58, RN59} gives the foundation for basic research in the field of epidemic modeling. These models were structured based on class age that means the time passed since an individual got infection. Katzmann and Dietz \cite{RN55} assume that maternal antibodies and immunity obtained from vaccination decay exponentially. Optimal age is determined for single vaccination during whole lifetime. In \cite{MR1814049} Hethcote considers many mathematical models of infectious diseases spread and also applied those models to some specific infectious diseases. Hadeler and Muller \cite{MR2319557} consider optimal vaccination and optimal harvesting problems together for age structured and size structured population models. More recently Chekroun and  Kuniya \cite{MR4039141} studied the global asymptotic behavior of an age structured SIR model with diffusion in a general n-dimensional bounded spatial domain under the homogeneous Dirichlet boundary condition. Without age structure, there are some works on global stability analysis of population models and analysis of epidemics models. \cite{syed1} studied global asymptotic stability of a dynamical system by constructing an appropriate Lyapunov function. \cite{syed2} analysed the effects of non pharmaceutical interventions for COVID-19 and \cite{syed3} developed a numerical scheme for fractional order SIR model by using Bernstein wavelets. \cite{syed4} considered chaos in an epidemic model and \cite{syed5} derived positive almost periodic solutions for population model with delay term.  There is enough literature available on age structured population models and optimal vaccination patterns, more details can be found in \cite{MR1647874,MR1971508, RN53 ,RN56,  MR4039141, MR4032666, MR4030068}.
		
Many times basic reproduction number can be calculated analytically by two methods. In one method, we need to find the threshold condition above which an endemic equilibrium exists and to interpret this condition in terms of basic reproduction number as $R_{0}>1$. In other method, local stability analysis of disease free equilibrium point is being done and the basic reproduction number can be evaluated from the threshold condition at which equilibrium point changes its asymptotic stability to instability. The $R_{0}$ obtained from both the methods will be same for SIR endemic models. Similar methods can be extended to obtain basic reproduction number in age structured SIR models. We consider an age structured SEIR model with quarantine of those individuals who are coming from the most affected areas (areas which have more number of cases ). We assume that disease also have some incubation period i.e. patient may not show symptoms just after getting infected. We assume that individuals in exposed class either enter into infected class or into quarantine (if they are traveling from areas which are most affected). Individuals from quarantine enter either into infected or recovered class. We derive steady state solutions to our model and analyzed the steady state solution with the assumption that force of infection involves separable mixing. Constructing an appropriate Liapunov function help us to check the stability of disease free equilibrium point. We also give vaccination to susceptible and exposed individuals and assume that there are no individuals in quarantine during vaccination programme. We also assume that due to strict government policies we are able to recognize those individuals who have already passed through the disease and acquired immunity. We also assume that infected individuals can be recognized and need not to be vaccinated. We show that individuals can be  vaccinated optimally at atmost three age classes.

Our paper is divided into $5$ sections, we formulate our model in section $2$ and defined various parameters used in the model.  Section $3$ is devoted to the study of  steady state solution to the given model. By constructing an appropriate Liapunov function, stability of disease free equilibrium is checked. We also find the average age of infection and basic reproduction number for constant parameters in this section. In section $4$, we formulate our new problem with vaccination and show that individuals can be optimally vaccinated at atmost three age classes. Last section is devoted to the overall discussion of our results.

\mysection{Model Formulation}
Let $U(a,t)$ be the age density of individuals of age $a$ at time $t$. $\mu(a)$ and $\beta(a)$ be age dependent mortality and fertility rates respectively.  Then the evolution of $U(a,t)$ can be described by the following McKendrick-Von Foerster PDE with boundary and initial conditions:

\begin{equation} \label{2.1}
\begin{cases}
  \frac{\partial U(a,t)}{\partial t} + \frac{\partial U(a,t) }{\partial a}=  - \mu (a)U(a,t) \quad (a,t) \in (0,\infty) \times (0,\infty)   \\
 U(0,t)= \int_{0}^{\infty}  \beta(a)U(a,t)da \quad t \in (0, \infty) \\
 U(a,0)=U_{0}(a) \quad a \in (0,\infty).\\
\end{cases}
\end{equation}
$U(0,t)$ is the number of newborns per unit time at time $t$. We assume that the mortality rate $ \mu \in L_{loc}^{1}([0,\infty))$ with $\int_{0}^{\infty} \mu(a) da= + \infty$ and the fertility rate $\beta \in L^{\infty}(0,\infty).$   $e^{-\int_{0}^{a} \mu(s) ds}$ is the proportion of individuals who are still living at age $a$ and $\int_{0}^{\infty} \beta(a)e^{-\int_{0}^{a} \mu(s) ds} da$ represents the net reproduction rate. Let us assume that the net reproduction rate is $1$. In this case steady state solution is given by \\ $U(a,t)=U(a)=\beta_{0} e^{-\int_{0}^{a} \mu(\tau) d \tau}$, where $\beta_{0}$ is given by $$ \beta_{0}= \frac{1}{\int_{0}^{\infty} e^{-\int_{0}^{a}\mu(\tau) d \tau} da}. $$
Let $S(a,t),E(a,t),Q(a,t), I(a,t)$ and $R(a,t)$ be the densities of susceptible, exposed,quarantined, infective and recovered individuals of age $a$ at time $t$. $k(a,b)$ is the age dependent transmission coefficient which describes the contact process between susceptible and infective individuals i.e. $k(a,b)S(a,t)I(b,t) da db$ is the number of susceptibles aged in $(a,a+da)$ that contract the disease by means of a suitable contact with an infective aged in $(b,b+db)$ . Let us assume that the force of infection i.e. the per capita rate of susceptibles to be infected at age $a$ at time $t$ is given in the following functional form $$ \phi(a,t) = \int_{0}^{\infty} k(a,\sigma) I(\sigma,t) d \sigma. $$ Then the spread of disease is described by the following system of partial differential equations
\begin{equation}
\begin{cases}
  \frac{\partial S(a,t)}{\partial t} + \frac{\partial S(a,t) }{\partial a}= -\phi(a,t)S(a,t) - \mu (a)S(a,t)    \\
  \frac{\partial E(a,t)}{\partial t} + \frac{\partial E(a,t) }{\partial a}= \phi(a,t)S(a,t)-\mu_{1}E(a,t)-q_{1}E(a,t) - \mu (a)E(a,t)    \\
    \frac{\partial Q(a,t)}{\partial t} + \frac{\partial Q(a,t) }{\partial a}= q_{1}E(a,t)-\gamma_{1}Q(a,t) - \gamma_{2} Q(a,t) -\mu(a)Q(a,t)    \\
  \frac{\partial I(a,t)}{\partial t} + \frac{\partial I(a,t) }{\partial a}= \mu_{1}E(a,t)+\gamma_{1} Q(a,t)  - \gamma I(a,t) -\mu(a)I(a,t)
  \\
  \frac{\partial R(a,t)}{\partial t} + \frac{\partial R(a,t) }{\partial a}= \gamma I(a,t)+\gamma_{2}Q(a,t)  - \mu (a)R(a,t)
  \\
S(0,t)= \int_{0}^{\infty}  \beta(a)(S(a,t)+E(a,t)+Q(a,t)+I(a,t)+R(a,t))da \\ E(0,t)=Q(0,t)=I(0,t)=R(0,t)=0 \\
 S(a,0)=S_{0}(a),E(a,0)=E_{0}(a),Q(a,0)=Q_{0}(a),I(a,0)=I_{0}(a)~\text{and}~ R(a,0)=R_{0}(a).    \\

\end{cases}
\end{equation}

\vspace{0.1cm}
\begin{center}
    \begin{tabular}{ | l | p{15cm} |}
    \hline
    $\mu(a)$ &  Age dependent natural mortality rate \\
    \hline
    $\mu_{1}$ & Progression rate from exposure to onset of symptoms \\
    \hline
    $q_{1}$  & Proportion of recruitment of asymptomatic persons for quarantine    \\
\hline
$\gamma_{1}$  & Proportion of individuals entering into infected class    \\
\hline
$\gamma_{2}$  & Proportion of individuals recovered after quarantine period    \\
\hline
$\gamma$  & Recovery rate of infected population    \\
\hline
    \end{tabular}
\end{center}
   \vspace{0.1cm}
   \begin{center}
 \includegraphics[width=15cm,height=9cm]{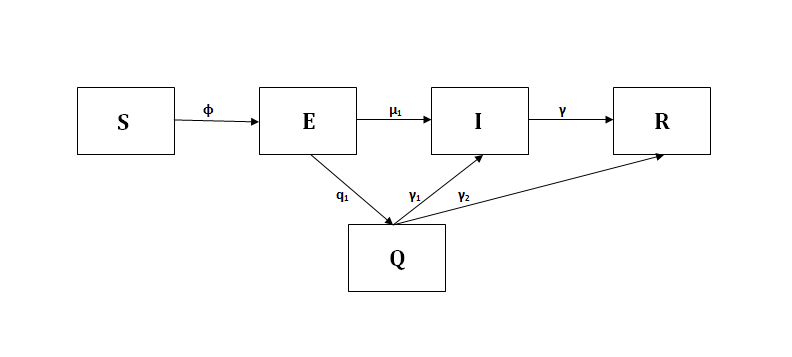}
\end{center}


 Let $s(a,t), e(a,t),q(a,t),i(a,t)$ and $r(a,t)$ be defined in the following way
 $$ s(a,t)=\frac{S(a,t)}{U(a,t)},e(a,t)=\frac{E(a,t)}{U(a,t)},q(a,t)=\frac{Q(a,t)}{U(a,t)}, i(a,t)=\frac{I(a,t)}{U(a,t)}~ \text{and}~ r(a,t)=\frac{R(a,t)}{U(a,t)} $$
 and the force of infection is given by
 $$ \phi(a,t) = \int_{0}^{\infty} r(a,\sigma) U(\sigma) i(\sigma,t) d \sigma. $$
 Then our new system becomes
 \begin{equation} \label{2.3}
\begin{cases}
  \frac{\partial s(a,t)}{\partial t} + \frac{\partial  s(a,t) }{\partial a}= -\phi(a,t) s(a,t)  \\
  \frac{\partial e(a,t)}{\partial t} + \frac{\partial e(a,t) }{\partial a}= \phi(a,t) s(a,t)-\mu_{1}e(a,t)-q_{1}e(a,t) \\
   \frac{\partial q(a,t)}{\partial t} + \frac{\partial q(a,t) }{\partial a}= q_{1}e(a,t)-\gamma_{1}q(a,t)-\gamma_{2}q(a,t) \\
   \frac{\partial i(a,t)}{\partial t} + \frac{\partial i(a,t) }{\partial a}= \mu_{1}e(a,t)+\gamma_{1}q(a,t)-\gamma i(a,t) \\
  \frac{\partial r(a,t)}{\partial t} + \frac{\partial r(a,t) }{\partial a}= \gamma i(a,t) + \gamma_{2}q(a,t)
  \\
  s(0,t)= 1, ~ e(0,t)=q(0,t)=0=i(0,t)=r(0,t)=0 \\
 s(a,0)=s_{0}(a),~e(a,0)=e_{0}(a),~q(0,a)=q_{0}(a),~i(a,0)=i_{0}(a)~\text{and}~ r(a,0)=r_{0}(a)
  \\
 s(a,t)+e(a,t)+q(a,t)+i(a,t)+r(a,t)=1. \\
\end{cases}
\end{equation}
 So, new transformations reduced our system into a simpler form i.e. boundary conditions now become constant and there is no term involving natural mortality rate.

\mysection{Steady state solutions}
Our system in steady state can be written as
\begin{equation} \label{1}
\begin{cases}
   \frac{d  s(a) }{da}= -\phi(a) s(a)    \\
  \frac{d  e(a) }{da}= \phi(a) s(a)-\mu_{1}e(a)-q_{1}e(a)    \\
  \frac{d  q(a) }{da}= q_{1} e(a)-\gamma_{1}q(a)-\gamma_{2}q(a)    \\
    \frac{d  i(a) }{da}= \mu_{1} e(a)+\gamma_{1}q(a)-\gamma i(a)    \\
      \frac{d  r(a) }{da}= \gamma i(a)+\gamma_{2}q(a)    \\
  s(0)= 1, e(0)=q(0)=i(0)=r(0)=0 \\
 s(a,t)+e(a,t)+q(a,t)+i(a,t)+r(a,t)=1 \\
  \phi(a) = \int_{0}^{\infty} k(a,\sigma) U(\sigma) i(\sigma) d \sigma .
\end{cases}
\end{equation}
Steady state solutions can be obtained as
\begin{eqnarray} \label{2}
s(a) &=& \exp \left( -\int_{0}^{a} \phi(\sigma) d \sigma \right)  \\
\label{3}
 e(a) &=& \int_{0}^{a} \phi(\sigma) e^{-\int_{0}^{\sigma} \phi(z) dz} e^{(\mu_{1}+q_{1})(\sigma -a)} d \sigma   \\
 \label{4}
 q(a)&=& \int_{0}^{a}q_{1}e(\sigma) e^{(\gamma_{1}+\gamma_{2})(\sigma-a)} d \sigma  \\
\label{5}
 i(a) &=& \int_{0}^{a} (\mu_{1} e(\sigma) + \gamma_{1}q(\sigma)) e^{\gamma (\sigma -a)} d \sigma.
 \end{eqnarray}
The force of infection is given by

\begin{eqnarray} \label{f1}
\phi(a) = \int_{0}^{\infty} k(a,\sigma) U(\sigma) i(\sigma) d \sigma 
 \end{eqnarray}
which depends on number of infected individuals and so explicitly depends on number of exposed and quarantined individuals. Using (\ref{3}),  (\ref{4}) can be written as
\begin{eqnarray} \label{6}
q(a) = \int_{0}^{a} \left[ q_{1} \left( \int_{0}^{\eta} \phi(y) e^{- \int_{0}^{y} \phi(z) dz} e^{(\mu_{1}+q_{1})(y-\eta)} dy \right) e^{(\gamma_{1}+\gamma_{2})(\eta -a)} \right] d \eta
\end{eqnarray}

and using (\ref{6}), (\ref{5}) can be written as \\

\begin{dmath} \label{7}
i(a) = \int_{0}^{a} \left[ \mu_{1} \left( \int_{0}^{\sigma} \phi(y) e^{-\int_{0}^{y} \phi(z) dz} e^{(\mu_{1}+q_{1})(y-\sigma)} dy \right)
+ \gamma_{1} \left\{ \int_{0}^{\sigma} q_{1} \left( \int_{0}^{\eta} \phi(y) e^{-\int_{0}^{y}\phi(z) dz} e^{(\mu_{1}+q_{1})(y-\eta)} dy \right) e^{(\gamma_{1}+\gamma_{2})(\eta-\sigma)} d \eta \right\} \right] e^{\gamma(\sigma-a)} d \sigma.
\end{dmath}
Therefore the force of infection is given by \\
\begin{dmath}
\phi(a) = \int_{0}^{\infty} k(a,b)U(b) \int_{0}^{b} \left[ \mu_{1} \left( \int_{0}^{\sigma} \phi(y) e^{-\int_{0}^{y} \phi(z) dz} e^{(\mu_{1}+q_{1})(y-\sigma)} dy \right)
+ \gamma_{1} \left\{ \int_{0}^{\sigma} q_{1} \left( \int_{0}^{\eta} \phi(y) e^{-\int_{0}^{y}\phi(z) dz} e^{(\mu_{1}+q_{1})(y-\eta)} dy \right) e^{(\gamma_{1}+\gamma_{2})(\eta-\sigma)} d \eta \right\} \right] e^{\gamma(\sigma-b)} d \sigma db.
\end{dmath}
Let $k(a,b)=k_{1}(a)k_{2}(b)$ i.e. we are assuming separable mixing and $\phi(a) = hk_{1}(a)$, then
\begin{dmath}
1= \int_{0}^{\infty} k_{2}(b)U(b) \int_{0}^{b} \left[ \mu_{1} \left( \int_{0}^{\sigma}  k_{1}(y) e^{-\int_{0}^{y} hk_{1}(z) dz} e^{(\mu_{1}+q_{1})(y-\sigma)} dy \right)
+ \gamma_{1} \left\{ \int_{0}^{\sigma} q_{1} \left( \int_{0}^{\eta} k_{1}(y) e^{-\int_{0}^{y}hk_{1}(z) dz} e^{(\mu_{1}+q_{1})(y-\eta)} dy \right) e^{(\gamma_{1}+\gamma_{2})(\eta-\sigma)} d \eta \right\} \right] e^{\gamma(\sigma-b)} d \sigma db.
\end{dmath}
 For $h=0$, define

 \begin{dmath} \label{10}
 R_{0} = \int_{0}^{\infty} k_{2}(b)U(b) \int_{0}^{b} \left[ \mu_{1} \left( \int_{0}^{\sigma}  k_{1}(y) e^{-\int_{0}^{y} hk_{1}(z) dz} e^{(\mu_{1}+q_{1})(y-\sigma)} dy \right)
+ \gamma_{1} \left\{ \int_{0}^{\sigma} q_{1} \left( \int_{0}^{\eta} k_{1}(y) e^{-\int_{0}^{y}hk_{1}(z) dz} e^{(\mu_{1}+q_{1})(y-\eta)} dy \right) e^{(\gamma_{1}+\gamma_{2})(\eta-\sigma)} d \eta \right\} \right] e^{\gamma(\sigma-b)} d \sigma db.
 \end{dmath}
 Our task is show that the stability of disease free equilibrium point depends   on $R_{0}$.
Now, let us construct the Liapunov function \\
\begin{dmath}
V = \int_{0}^{\infty} \left[ \alpha_{1}(a)e(a,t) +\alpha_{2}(a)q(a,t)+\alpha_{3}(a)i(a,t) \right] da
\end{dmath}

 \begin{dmath*}
 \dot{V} = \int_{0}^{\infty} \left[ \alpha_{1}(a) \left( \phi(a,t)s-\mu_{1}e-q_{1}e-\frac{\partial e}{\partial a} \right)+\alpha_{2}(a) \left( q_{1}e - \gamma_{1}q - \gamma_{2}q - \frac{\partial q}{\partial a} \right)+\alpha_{3}(a) \left( \mu_{1}e + \gamma_{1}q - \gamma i - \frac{\partial i}{\partial a} \right) \right]
 \end{dmath*}
 \begin{dmath*}
 = \int_{0}^{\infty} \left[ \alpha_{1}(a) \phi(a,t)s + e ( \alpha'_{1}(a)-\mu_{1}(a) \alpha_{1}(a) - q_{1} \alpha_{1}(a) + q_{1} \alpha_{2}(a)+\mu_{1} \alpha_{3}(a)) + q ( \alpha'_{2}(a) - \gamma_{1}\alpha_{2}(a) - \gamma_{2} \alpha_{2}(a) + \gamma_{1} \alpha_{3}(a) ) + i ( \alpha'_{3}(a)- \gamma \alpha_{3}(a)) \right] da
 \end{dmath*}
Choose $\alpha_{1}, \alpha_{2}$ such that  coefficients of $e$ and $q$ are zero i.e.

\begin{eqnarray}
 \alpha'_{1}(a)-\mu_{1}(a) \alpha_{1}(a) - q_{1} \alpha_{1}(a) + q_{1} \alpha_{2}(a)+\mu_{1} \alpha_{3}(a) =0 \\
  \alpha'_{2}(a) - \gamma_{1}\alpha_{2}(a) - \gamma_{2} \alpha_{2}(a) + \gamma_{1} \alpha_{3}(a)=0
\end{eqnarray}
Let us choose initial conditions as
\begin{eqnarray}
\alpha_{1}(0) &=& \int_{0}^{\infty} (q_{1} \alpha_{2}(z)+\mu_{1} \alpha_{3}(z)) e^{-(\mu_{1}+q_{1})z} dz \\
\alpha_{2}(0) &=& \int_{0}^{\infty} \gamma_{1} \alpha_{3}(z) e^{-(\gamma_{1}+\gamma_{2})z} dz.
\end{eqnarray}
 Then the solutions are given by
 \begin{eqnarray} \label{16}
 \alpha_{1}(a) &=& \int_{a}^{\infty} (q_{1} \alpha_{2}(z)+\mu_{1} \alpha_{3}(z)) e^{(\mu_{1}+q_{1})(a-z)} dz \\
 \label{17} \alpha_{2}(a) &=& \int_{a}^{\infty} \gamma_{1} \alpha_{3}(z) e^{(\gamma_{1}+\gamma_{2})(a-z)} dz.
 \end{eqnarray}
 Using (\ref{17}) in (\ref{16}), we get
 \begin{eqnarray}
 \alpha_{1}(a) =  \int_{a}^{\infty} \left[ q_{1} \int_{z}^{\infty} \gamma_{1} \alpha_{3}(y) e^{(\gamma_{1}+\gamma_{2})(z-y)} dy+\mu_{1} \alpha_{3}(z) \right] e^{(\mu_{1}+q_{1})(a-z)} dz
 \end{eqnarray}
 Therefore,
 \begin{dmath}
 \dot{V} = \int_{0}^{\infty} s k_{1}(a) \int_{0}^{\infty} k_{2}(b) i(b,t) U(b)db
 \int_{a}^{\infty} \left[ q_{1} \int_{z}^{\infty} \gamma_{1} \alpha_{3}(y) e^{(\gamma_{1}+\gamma_{2})(z-y)} dy+\mu_{1} \alpha_{3}(z) \right] e^{(\mu_{1}+q_{1})(a-z)} dz da + \int_{0}^{\infty} i \left( \alpha'_{3}(a)-\gamma \alpha_{3}(a) \right) da.
 \end{dmath}
 Choose $\alpha_{3}$ such that
 \begin{dmath*}
 \int_{0}^{\infty} i( \alpha'_{3}(a) - \gamma \alpha_{3}(a) ) da = -\int_{0}^{\infty}i k_{1}(a) U(a) da
 \end{dmath*}
 $$ \text{i.e}~ \alpha'_{3}(a) - \gamma \alpha_{3}(a) = -k_{2}(a) U(a). $$
 Choose initial condition as
 $$ \alpha_{3}(0) = \int_{0}^{\infty} k_{2}(x) U(x) e^{- \gamma x} dx .$$
 Then $\alpha_{3}$ can be evaluated as
 \begin{dmath}
 \alpha_{3}(a) = \int_{a}^{\infty} k_{2}(x) U(x) e^{\gamma(a-x)} dx
 \end{dmath}
 So, derivative of Liapunov function can be written as
\begin{dmath} \label{21}
\dot{V} = \left[ \int_{0}^{\infty} s k_{1}(a) \int_{a}^{\infty} \left \{ q_{1} \int_{z}^{\infty} e^{\gamma y} \left( \int_{y}^{\infty}k_{2}(x)U(x) e^{-\gamma x} dx \right) \gamma_{1} e^{(\gamma_{1}+\gamma_{2})(z-y)}dy + \mu_{1} e^{\gamma z} \int_{z}^{\infty} k_{2}(x)U(x) e^{-\gamma x} dx \right \} e^{(\mu_{1}+q_{1})(a-z)} dz da -1  \right] \int_{0}^{\infty} k_{1}(b)i(b,t)U(b) db
 = \left[ \int_{0}^{\infty} s k_{1}(a) e^{(\mu_{1}+q_{1})a} \int_{a}^{\infty} q_{1}e^{(\gamma_{1}+\gamma_{2})z} e^{-(\mu_{1}+q_{1})z} \int_{z}^{\infty} \gamma_{1}e^{\gamma y} e^{-(\gamma_{1}+\gamma_{2})y} \int_{y}^{\infty} k_{2}(x)U(x)e^{-\gamma x} dx dy dz da + \int_{0}^{\infty} k_{1}(a)e^{(\mu_{1}+q_{1})a} \int_{a}^{\infty} \mu_{1} e^{\gamma z} e^{-(\mu_{1}+q_{1})z} \int_{z}^{\infty} k_{2}(x)U(x) e^{- \gamma x} dx dz da - 1 \right] \int_{0}^{\infty} k_{1}(b)i(b,t)U(b) db
 = \left[ \mathcal{I} - 1 \right] \int_{0}^{\infty} k_{1}(b)i(b,t)U(b) db
\end{dmath}
where
\begin{dmath} \label{n3.22}
\mathcal{I} = \int_{0}^{\infty} s k_{1}(a) e^{(\mu_{1}+q_{1})a} \int_{a}^{\infty} q_{1}e^{(\gamma_{1}+\gamma_{2})z} e^{-(\mu_{1}+q_{1})z} \int_{z}^{\infty} \gamma_{1}e^{\gamma y} e^{-(\gamma_{1}+\gamma_{2})y} \int_{y}^{\infty} k_{2}(x)U(x)e^{-\gamma x} dx dy dz da + \int_{0}^{\infty} k_{1}(a)e^{(\mu_{1}+q_{1})a} \int_{a}^{\infty} \mu_{1} e^{\gamma z} e^{-(\mu_{1}+q_{1})z} \int_{z}^{\infty} k_{2}(x)U(x) e^{- \gamma x} dx dz da
= \mathcal{I}_{1} + \mathcal{I}_{2}.
\end{dmath}
The region of integration for $\mathcal{I}_{1}$ is
\begin{dmath}
\{ (x,y,z,a) \in \mathbb{R}^{4} ~|~ 0 \le a \le z \le y \le x \}
\end{dmath}
and for $\mathcal{I}_{2}$ is
\begin{dmath}
\{ (x,z,a) \in \mathbb{R}^{3} ~|~ 0 \le a \le z  \le x \}.
\end{dmath}
 Observe that
 \begin{dmath} \label{R_{0}3.22}
 R_{0} = \int_{0}^{\infty} k_{2}(b) U(b) e^{-\gamma b} \int_{0}^{b} \mu_{1} e^{\gamma \sigma} e^{-(\mu_{1}+q_{1})\sigma} \int_{0}^{\sigma} k_{1}(y)e^{(\mu_{1}+q_{1})y} dy d \sigma db + \int_{0}^{\infty}k_{2}(b)U(b) e^{-\gamma b} \int_{0}^{b} \gamma_{1} e^{\gamma \sigma} e^{-(\gamma_{1}+\gamma_{2})\sigma} \int_{0}^{\sigma}q_{1} e^{(\gamma_{1}+\gamma_{2})\eta} e^{-(\mu_{1}+q_{1})\eta} \int_{0}^{\eta} k_{1}(y) e^{(\mu_{1}+q_{1})y} dy d \eta d \sigma db
 \end{dmath}
 Also after changing the role of dummy variables, we get
 \begin{dmath} \label{nR_{0}3.22}
 R_{0} = \int_{0}^{\infty}k_{2}(x)U(x) e^{-\gamma x} \int_{0}^{x} \gamma_{1} e^{\gamma y} e^{-(\gamma_{1}+\gamma_{2})y} \int_{0}^{y}q_{1} e^{(\gamma_{1}+\gamma_{2})z} e^{-(\mu_{1}+q_{1})z} \int_{0}^{z} k_{1}(a) e^{(\mu_{1}+q_{1})a} da dz dy dx  + \int_{0}^{\infty} k_{2}(x) U(x) e^{-\gamma x} \int_{0}^{x} \mu_{1} e^{\gamma z} e^{-(\mu_{1}+q_{1})z} \int_{0}^{z} k_{1}(a)e^{(\mu_{1}+q_{1})a} da dz dx
 = \mathcal{R}_{1} + \mathcal{R}_{2} .
 \end{dmath}
 The region of integration for $\mathcal{R}_{1}$ is
\begin{dmath}
\{ (x,y,z,a) \in \mathbb{R}^{4} ~|~ 0 \le a \le z \le y \le x \}
\end{dmath}
and for $\mathcal{R}_{2}$ is
\begin{dmath}
\{ (x,z,a) \in \mathbb{R}^{3} ~|~ 0 \le a \le z  \le x \}.
\end{dmath}

 Since $s \le 1$, after changing the order of integration, the integral in (\ref{n3.22}) with $s=1$ equal to $R_{0}$ defined in (\ref{nR_{0}3.22}). Therefore,
 \begin{eqnarray} \label{3.23}
 \dot{V} \le (R_{0}-1)\int_{0}^{\infty} k_{1}(b)i(b,t)U(b) db \le 0 ~~~ \text{if}~ R_{0} \le 1.
 \end{eqnarray}
	
Therefore, through the level sets of Liapunov function $V$, solutions of (\ref{2.3})	move downward as long as they do not stop on the set where $\dot{V} =0$. From (\ref{3.23}) it is clear that the set $\dot{V} = 0$ is the boundary of the feasible region with $i=0$, but $$\frac{di}{dt}(a(t),t)  = \mu_{1}e+\gamma_{1}q$$ on this boundary. So, number of infected individuals $i$ moves off this boundary unless $e=q=0$. If $e=q=i=0$, there are no exposed or quarantined or infected individuals, then there would be no recovered individuals. So, every individual will be susceptible and disease free steady state is the only  positively invariant subset of the set with $\dot{V}=0$. Similarly, if we have $R_{0}>1$, then for the points near to the disease free equilibrium point the derivative of Liapunov function will be positive with $s$ close to $1$ and for some age the number of infected individuals will be non zero. So, disease free equilibrium will be unstable in this case. For finite maximum age, we can apply the same analysis as used by \cite{MR1814049} to show that all the paths in the feasible region will approach the disease free equilibrium point for the case $R_{0} <1$. And also can show that the disease free equilibrium point is unstable for the case $R_{0} >1.$
\subsection{Average Age of Infection}
Here, we will find the expression for the average age of infection i.e. the average time spent in the susceptible group before becoming infected. The steady state age distribution of the population is given by 
$$ U(a) = \beta_{0} e^{-\int_{0}^{a} \mu(\tau) d \tau}, ~\text{where}~  \beta_{0}= \frac{1}{\int_{0}^{\infty} e^{-\int_{0}^{a}\mu(\tau) d \tau} da}. $$ 
The age distribution for a particular birth cohort will be given by 
$$ \frac{e^{-\int_{0}^{a} \mu(z)dz}}{\int_{0}^{\infty}e^{-\int_{0}^{a}\mu(z)dz}da}. $$ 
Thus the leaving rate of individuals in a particular birth cohort from susceptible class due to infection is given by 
$$ \phi(a)s(a) \frac{e^{-\int_{0}^{a} \mu(z)dz}}{\int_{0}^{\infty}e^{-\int_{0}^{a}\mu(z)dz}da} $$
where $\phi(a)$ and $s(a)$ are given by (\ref{f1}) and (\ref{2}) respectively. Since force of infection depends on the number of infected individuals and the number of infected individuals depend on number of quarantine individuals. So this rate also depend on quarantine individuals which is an important factor to contain the spread of disease. Therefore, average age of infection is given by 
\begin{equation}
\mathcal{A} = \frac{\int_{0}^{\infty}a \phi(a) s(a) e^{-\int_{0}^{a} \mu(z)dz} da}{\int_{0}^{\infty} \phi(a)s(a) e^{-\int_{0}^{a} \mu(z) dz}da}.
\end{equation} 

\subsection{Negative Exponential Survival}
In this section, let us assume that mortality rate $\mu(a)$ is independent of age $a$, then the steady state age distribution is given by 
$$ U(a) = \mu e^{-\mu a}. $$ Also let us assume that contact rate does not depends on the ages of susceptible and infective individuals, so let $k_{1}(a)=1$ and $k_{2}(a)=k_{2}$. With these assumptions $\phi(a)$ will be constant and let this constant be denoted by $\phi$. Now the basic reproduction number will be given by  
 \begin{dmath} \label{CR_{0}3.22}
 R_{0} = \int_{0}^{\infty} k_{2} \mu e^{-\mu a} e^{-\gamma b} \int_{0}^{b} \mu_{1} e^{\gamma \sigma} e^{-(\mu_{1}+q_{1})\sigma} \int_{0}^{\sigma} e^{(\mu_{1}+q_{1})y} dy d \sigma db + \int_{0}^{\infty}k_{2} \mu e^{-\mu b} e^{-\gamma b} \int_{0}^{b} \gamma_{1} e^{\gamma \sigma} e^{-(\gamma_{1}+\gamma_{2})\sigma} \int_{0}^{\sigma}q_{1} e^{(\gamma_{1}+\gamma_{2})\eta} e^{-(\mu_{1}+q_{1})\eta} \int_{0}^{\eta}  e^{(\mu_{1}+q_{1})y} dy d \eta d \sigma db
 \end{dmath}
Let us denote the first integral in (\ref{CR_{0}3.22}) by $R_{1}$ and second by $R_{2}$. Then 
$$ R_{1} = \frac{\mu_{1}}{\mu_{1}+q_{1}} \left[ \frac{k_{2}\beta_{0}}{\gamma} \left\{ \frac{1}{\mu} -\frac{1}{\mu+\gamma} \right\}  - \frac{k_{2}\beta_{0}}{\gamma - \mu_{1}-q_{1}} \left\{ \frac{1}{\mu + \mu_{1}+q_{1}} -\frac{1}{\mu+\gamma} \right\}\right] $$
$$ = \frac{k_{2}\beta_{0} \mu_{1}}{(\mu_{1}+q_{1})(\mu+\gamma)} \left\{ \frac{1}{\mu} - \frac{1}{\mu+\mu_{1}+q_{1}} \right\}  $$
$$ = \frac{k_{2}\beta_{0} \mu_{1}}{\mu (\mu+ \gamma)(\mu + \mu_{1}+q_{1})}. $$
Because $\beta_{0}= \frac{1}{\int_{0}^{\infty} e^{-\int_{0}^{a}\mu(\tau) d \tau} da} = \mu$, we have 
$$R_{1} =\frac{k_{2} \mu_{1}}{ (\mu+ \gamma)(\mu + \mu_{1}+q_{1})}.  $$
Also, 
\begin{dmath}
R_{2} = \frac{1}{\mu_{1}+q_{1}} \left[ \frac{k_{2}\beta_{0} \gamma_{1}q_{1}}{\gamma(\gamma_{1}+\gamma_{2})} \left\{ \frac{1}{\mu} -\frac{1}{\mu + \gamma} \right\} -\frac{k_{2}\beta_{0}\gamma_{1}q_{1}}{(\gamma_{1}+\gamma_{2})(\gamma - \gamma_{1}-\gamma_{2})} \left\{ \frac{1}{\mu+\gamma_{1}+\gamma_{2}}- \frac{1}{\mu+\gamma} \right\}  -\frac{k_{2} \beta_{0}\gamma_{1}q_{1}}{(\gamma_{1}+\gamma_{2}-\mu_{1}-q_{1})(\gamma -\mu_{1}-q_{1})} \left\{ \frac{1}{\mu+\mu_{1}+q_{1}} -\frac{1}{\mu+\gamma} \right\} + \frac{k_{2}\beta_{0}\gamma_{1}q_{1}}{(\gamma_{1}+\gamma_{2}-\mu_{1}-q_{1})(\gamma-\gamma_{1}-\gamma_{2})} \left\{ \frac{1}{\mu+\gamma_{1}+\gamma_{2}} - \frac{1}{\mu+\gamma}  \right\}\right]
 = \frac{1}{\mu_{1}+q_{1}} \left[ \frac{k_{2}\beta_{0}\gamma_{1}q_{1}}{(\mu+\gamma)(\mu+\gamma_{1}+\gamma_{2})} \left\{  \frac{1}{\mu} - \frac{1}{\mu+\mu_{1}+q_{1}} \right\} \right]
 = \frac{k_{2} \gamma_{1} q_{1}}{(\mu+\gamma)(\mu+\gamma_{1}+\gamma_{2})(\mu+\mu_{1}+q_{1})}
\end{dmath}
where we have used the fact that $\beta_{0} = \mu.$
Therefore,
$$ R_{0} = \frac{k_{2} \mu_{1}}{ (\mu+ \gamma)(\mu + \mu_{1}+q_{1})} + \frac{k_{2} \gamma_{1} q_{1}}{(\mu+\gamma)(\mu+\gamma_{1}+\gamma_{2})(\mu+\mu_{1}+q_{1})} $$
$$ = \frac{k_{2}\gamma_{1}(q_{1}+\gamma_{1})}{(\mu+\gamma)(\mu+\mu_{1}+q_{1})(\mu+\gamma_{1}+\gamma_{2})} + \frac{k_{2}\mu_{1}}{(\mu+\gamma)(\mu+\mu_{1}+q_{1})}.  $$
This $R_{0}$ has same interpretation as we have for any SEQIR model without age structure.
Note that if there are no quarantine individuals, then $\gamma_{1}=\gamma_{2}=q_{1}=0$. So, basic reproduction number becomes 
$$R_{0} =\frac{k_{2}\mu_{1}}{(\mu+\gamma)(\mu+\mu_{1})}. $$
Above basic reproduction number is same as the  reproduction derived for SEIR model without age structure in \cite{MR1814049}. 
Let us find the average age of infection for this SEQIR model. In this section, steady state age distribution will be given by 
$$ U(a) = \mu e^{-\mu a}. $$ Also the age distribution for a particular birth cohort is given by $\mu e^{-\mu a}$. Thus the leaving rate of individuals in a particular birth cohort from susceptible class due to infection is given by 
$$ \phi s(a) \mu e^{-\mu a}. $$
Therefore, the average age of infection is given by 
$$\mathcal{A} = \frac{\int_{0}^{\infty}a \phi e^{-\phi a} e^{-\mu a} da}{\int_{0}^{\infty} \phi e^{-\phi a} e^{-\mu a} da}  $$
$$ = \frac{1}{\phi+\mu} .$$
This average age of infection matches with the average age of infection given for MSEIR model \cite{MR1814049} without passively immune infants.
\mysection{Optimal Vaccination}
Let $v(a)$ be the vaccination policy for age $a$ individuals and also assume that there are no  quarantined individuals during vaccination programme. We are assuming that due to government policies, it is possible to recognize those individuals who have already passed through the disease and acquired immunity. We are also not vaccinating the infected individuals. $\mathcal{V}(a)$ is the density of those age $a$ individuals who became immune after vaccination programme i.e. those individuals who are coming from susceptible and exposed class after vaccination policy.
Then the steady state equations can be written as  \\
\begin{equation} \label{1}
\begin{cases}
   \frac{d  S(a) }{da}= -hk_{1}(a) S(a)-v(a)S(a) -\mu(a)S(a)     \\
  \frac{d  E(a) }{da}= hk_{1}(a) S(a)-\mu_{1}E(a)-v(a)E(a)- \mu(a) E(a)    \\
    \frac{d  I(a) }{da}= \mu_{1} E(a) - \gamma I(a) - \mu(a) I(a)   \\
    \frac{d  R(a) }{da}= \gamma I(a) - \mu R(a) \\
      \frac{d  \mathcal{V}(a) }{da}= v(a)S(a)+v(a)E(a)-\mu(a)\mathcal{V}(a)    \\
  S(0)= \int_{0}^{\infty} \beta(a)(S(a)+E(a)+I(a)+R(a)+\mathcal{V}(a)) da \\
   E(0)=I(0)=R(0)=\mathcal{V}(0)=0 \\
   h = \int_{0}^{\infty} k_{2}(\sigma) I(\sigma) d \sigma .
\end{cases}
\end{equation}	
\\
\\

Steady state solutions are given by \\
\begin{eqnarray}
S(a) &=& U(a)D(a) \exp \left(-\int_{0}^{a}hk_{1}(z) dz \right) \\
E(a) &=& U(a)D(a) \int_{0}^{a} h U(\sigma) k_{1}(\sigma) \exp \left(-\int_{0}^{\sigma} hk_{1}(z) dz \right) e^{\mu_{1}(\sigma - a) }d \sigma \\
I(a) &=& U(a) \int_{0}^{\sigma} \mu_{1} \left \{ U(\sigma) D(\sigma) \int_{0}^{\sigma} h U(\eta)k_{1}(\eta) \exp \left( -\int_{0}^{\eta}hk_{1}(z) dz \right) e^{\mu_{1}(\eta - \sigma)} d \eta \right \} e^{\gamma(\sigma - a)} d \sigma \\
R(a) &=& U(a) \int_{0}^{a} \gamma I(\sigma) d \sigma
\end{eqnarray}
where $D(a) = \exp \left( -\int_{0}^{a}v(z) dz \right) = 1 - \int_{0}^{a} v(z) \exp \left(-\int_{z}^{a} v(\tau) d\tau \right) dz$. This expression for $D(a)$ is going to be very useful for further analysis. \\
Let us define the cost function by $\tilde{C}(v)$ with weight functions $g_{1}(a)$ and $g_{2}(a)$ which depend on the cost of moving individuals from susceptible and exposed class to immune class respectively.
\begin{eqnarray}
\tilde{C}(v) = \int_{0}^{\infty}v(a) (g_{1}(a)S(a)+g_{2}(a)E(a)) da.
\end{eqnarray}
Here, instead of considering basic reproduction number, we will consider number of infected individuals with weight function $f(a)$ which depends on the social impact of one diseased case. Therefore we will consider the weighted prevalence
\begin{eqnarray}
\tilde{F}(v) = \int_{0}^{\infty} f(a)I(a) da.
\end{eqnarray}
Also, let us define
\begin{eqnarray}
\tilde{H}(v) = \int_{0}^{\infty}k_{2}(a) I(a) da.
\end{eqnarray}
Then we can write
\begin{dmath*}
\tilde{C}(v) = \int_{0}^{\infty} v(a) \left[  g_{1}(a)U(a)D(a) \exp \left(-\int_{0}^{a} hk_{1}(z) dz \right) + g_{2}(a)U(a)D(a) \int_{0}^{a} hU(\sigma)k_{1}(\sigma) \exp \left( -\int_{0}^{\sigma}hk_{1}(z) dz \right) e^{\mu_{1}(\sigma -a)} d \sigma \right] da
\\
\tilde{F}(v) = \int_{0}^{\infty} \mu_{1}U(\sigma) D(\sigma) e^{(\gamma -\mu_{1}) \sigma} d \sigma \int_{\sigma}^{\infty} f(a) U(a) e^{-\gamma a} da \int_{0}^{\sigma} h U(\eta) k_{1}(\eta) \exp \left(-\int_{0}^{\eta} hk_{1}(z)\right) e^{\mu_{1}\eta} d \eta  \\
\tilde{H}(v) =  \int_{0}^{\infty} \mu_{1}U(\sigma) D(\sigma) e^{(\gamma -\mu_{1}) \sigma} d \sigma \int_{\sigma}^{\infty} k_{2}(a) U(a) e^{-\gamma a} da \int_{0}^{\sigma} h U(\eta) k_{1}(\eta) \exp \left(-\int_{0}^{\eta} hk_{1}(z)\right) e^{\mu_{1}\eta} d \eta.
	\end{dmath*}
	Let us define the kernels
	\begin{eqnarray*}
	C_{1}(a) &=& U(a) \left[ g_{1}(a)\exp \left(-\int_{0}^{a} hk_{1}(z) dz \right) + g_{2}(a) \int_{0}^{a}h U(\sigma)k_{1}(\sigma) \exp \left(-\int_{0}^{\sigma}hk_{1}(z) dz \right) e^{\mu_{1}(\sigma -a)} d \sigma \right] \\
	F_{1}(a) &=& \int_{0}^{\infty} \mu_{1} U(\sigma)e^{\gamma \sigma} d \sigma \int_{\sigma}^{\infty} f(s)U(s) e^{-\gamma s} ds \int_{0}^{\sigma} h U(\eta) k_{1}(\eta) \exp \left( -\int_{0}^{\eta} hk_{1}(z) dz \right) \\
	H_{1}(a) &=& \int_{0}^{\infty} \mu_{1} U(\sigma)e^{\gamma \sigma} d \sigma \int_{\sigma}^{\infty} K_{2}(s)U(s) e^{-\gamma s} ds \int_{0}^{\sigma} h U(\eta) k_{1}(\eta) \exp \left( -\int_{0}^{\eta} hk_{1}(z) dz \right)
	\end{eqnarray*}
In terms of $v$, our optimization problem will be highly nonlinear so, let us use the substitution $$ \psi(a) = \exp \left( -\int_{0}^{\infty}v(s) ds \right) v(a). $$
Then
$$ \int_{0}^{a} \psi(z) dz = 1- e^{-\int_{0}^{a}v(s) ds} $$ and the inverse transformation will be given by
$$ v(a) = \frac{\psi(a)}{1-\int_{0}^{a} \psi(a) ds} .$$
So, in terms of $\psi$, we have
\begin{eqnarray}
\tilde{C}(v) = C(\psi) = \int_{0}^{\infty} \psi(a) C_{1}(a) da \\
\tilde{F}(v) = \tilde{F}(0) - F(\psi),~F(\psi) = \int_{0}^{\infty}\psi(a)F_{1}(a)da \\
\tilde{H}(v) = \tilde{H}(0) - H(\psi),~H(\psi) = \int_{0}^{\infty}\psi(a)H_{1}(a)da.
\end{eqnarray}
Therefore, our optimal vaccination problem becomes \\
	$$ \text{Minimize} ~C(\psi) $$ subject to the conditions \\
	$$ \psi(a) \ge 0 ~\text{for}~0 \le a < \infty  $$
	$$ F(\psi) \ge \tilde{F}(0)- \bar{F} $$
	$$Q(\psi) = \int_{0}^{\infty} \psi(a)da \le 1 $$
	$$H(\psi) = \tilde{H}(0) -h $$
	where $\bar{F}$ is the upper bound for the weighted prevalence i.e. $\tilde{F}(v) \le \bar{F}$. \\
	Let us define the Lagrange functional \\
	$$ L(\psi, \lambda_{1},\lambda_{2},\lambda_{3}) = C(\psi) - \lambda_{1}(F(\psi)-\tilde{F}(0)+\bar{F}) - \lambda_{2}(Q(\psi)-1)-\lambda_{3}(H(\psi)-\tilde{H}(0)+h) $$
	where $\lambda_{1},\lambda_{2},\lambda_{3} \in \mathbb{R}_{+}.$ \\
	The Kuhn-Tucker conditions are given by \\
	\begin{eqnarray}
	 \psi(a) \ge 0, \\
 F(\psi) \ge \tilde{F}(0)-\bar{F}, \\
	 Q(\psi) \le 1, \\
	H(\psi) = \tilde{H}(0)-h, \\
	 \lambda_{1}(F(\psi)-\tilde{F}(0)+\bar{F}) = 0, \\
	\lambda_{2}(Q(\psi)-1)= 0,\\
	 \lambda_{3}(H(\psi)-\tilde{H}(0)+h) = 0, \\
	 C_{1}(a)-\lambda_{1}F_{1}(a) - \lambda_{2} - \lambda_{3}H_{1}(a) \le 0, \\
	 C(\psi)-\lambda_{1}F(\psi)-\lambda_{2}Q(\psi)-\lambda_{3}H(\psi)=0.
	\end{eqnarray}
	\textbf{Case 1:} $Q(\psi) <1$, then $\lambda_{2} = 0.$
	Firstly assume that $\lambda_{3} =0$, then we have
	$$ C_{1}(a)-\lambda_{1}F_{1}(a) \le 0, $$
	$$ C(\psi)-\lambda_{1}F(\psi)=0. $$
	So, we have
	$$\frac{C_{1}(a)}{F_{1}(a)} \le \lambda_{1} = \frac{C(\psi)}{F(\psi)}  =  \frac{\int_{0}^{\infty}C_{1}(a)\psi(a)da}{\int_{0}^{\infty}F_{1}(a)\psi(a)da}.$$
	So, $\lambda_{1}$ is equal to weighted arithmetic mean of quotients $\frac{C_{1}(a)}{F_{1}(a)}$ and by the virtue of lower bound, we can say that supremum $\lambda_{1}$ is assumed for atleast one age $a=A~$ i.e. \\
	$$ \lambda_{1} = \frac{C_{1}(A)}{F_{1}(A)}. $$ So, at age $A$ vaccination policy can be applied i.e optimal policy has the form of a delta peak \\
	$$ v(a) = c_{1} \delta_{A}(a). $$
The coefficient $c_{1}$ will give the intensity of the vaccination policy. Even if maximum value is attained at multiple points, then convex combination of the optimal policies is optimal. \\
Now, let us consider the case when $\lambda_{3}$ may not be zero, then \\
$$ \frac{C_{1}(a)-\lambda_{1}F_{1}(a)}{H_{1}(a)} \le \lambda_{3} = \frac{C(\psi)-\lambda_{1}F(\psi)}{H(\psi)}  = \frac{\int_{0}^{\infty} (C_{1}(a)-\lambda_{1}F_{1}(a))\psi(a)da}{\int_{0}^{\infty}H_{1}(a)\psi(a)da}. $$
Using the same argument, we can show the existence of optimal age policy $A$ such that \\
$$ \lambda_{3} = \frac{C_{1}(A)-\lambda_{1}F_{1}(A)}{H_{1}(A)}. $$
Here, in this case there are choices for $\lambda_{1}$, there are two such ages $A_{1}$ and $A_{2}$ and the optimal policy $v$ is a combination of two delta peaks at ages $A_{1}$ and $A_{2}$. This $v$ is a candidate for an optimal two age policy. So, we can have either two age or one age optimal policy. \\
\textbf{Case 2:} $Q(\psi)=1$, then $\lambda_{2}$ may not vanish. Firstly assume that $\lambda_{3}=0$, then we have
$$ C_{1}(a)-\lambda_{1}F_{1}(a) \le \lambda_{2} = C(\psi)-\lambda_{1}F(\psi) .$$
So, there exist $a=A$ such that $$\lambda_{2} = C_{1}(A)-\lambda_{1}F_{1}(A). $$ Then in similar manner either we have one age or two age policy. \\
 Now, let us again consider the case when $\lambda_{3}$ may not be zero, then we have
 $$\frac{C_{1}(a)-\lambda_{1}F_{1}(a)-\lambda_{2}}{H_{1}(a)} \le \lambda_{3} = \frac{C(\psi)-\lambda_{1}F(\psi)-\lambda_{2}}{H(\psi)}. $$
 So, in this case we can have either one age or two age or three age policy. \\
 If $v$ is a one age policy then we can write \\
 $$ v(a) = c_{1} \delta_{A}(a)$$
 where $A$ is the age at which vaccination is applied and $c_{1}$ is the intensity of vaccination. Now, \\
 \[
\int_{0}^{a} v(s)ds =
  \begin{cases}
     0 & a < A \\
      c_{1} & a > A.
   \end{cases}
\]

 \[
e^{-\int_{0}^{a} v(s)ds} =
  \begin{cases}
     1 & a < A \\
      e^{-c_{1}} & a > A.
   \end{cases}
\]
  Therefore, $$ \psi(a) = - \frac{d}{da} \left(e^{-\int_{0}^{a}v(s)ds}\right) = (1-e^{-c_{1}}) \delta_{A}(a). $$
So, $$C(\psi) = (1-e^{-c_{1}})C_{1}(A) $$
$$F(\psi) = (1-e^{-c_{1}})F_{1}(A)  $$
$$ H(\psi) = (1-e^{-c_{1}})H_{1}(A).  $$	
These expressions can help us to reduce our optimization problem into a finite dimensional optimization problem. \\
If $$v(a) = c_{1}\delta_{A_{1}}(a)+c_{2}\delta_{A_{2}}(a) $$
then in similar manner
$$ \psi(a) = (1-e^{-c_{1}})\delta_{A_{1}}(a)+(e^{-c_{1}}-e^{-c_{1}-c_{2}}) \delta_{A_{2}}(a). $$
Two age policy will satisfy either $Q(\psi)<1$ or $Q(\psi)=1$. If $Q(\psi) <1$, then
$$ \psi(a) = (1-e^{-c_{1}})\delta_{A_{1}}(a)+e^{-c_{1}}(1-e^{-c_{2}})\delta_{A_{2}}(a) $$
$$C(\psi) = (1-e^{-c_{1}})C_{1}(A_{1})+e^{-c_{1}}(1-e^{-c_{2}})C_{1}(A_{2})  $$
$$F(\psi) = (1-e^{-c_{1}})F_{1}(A_{1})+e^{-c_{1}}(1-e^{-c_{2}})F_{1}(A_{2})  $$
$$H(\psi) = (1-e^{-c_{1}})H_{1}(A_{1})+e^{-c_{1}}(1-e^{-c_{2}})H_{1}(A_{2}).  $$
If $Q(\psi) =1$, then $c_{2} = \infty$ and
$$ \psi(a) = (1-e^{-c_{1}})\delta_{A_{1}}(a)+e^{-c_{1}}\delta_{A_{2}}(a) $$
$$C(\psi) = (1-e^{-c_{1}})C_{1}(A_{1})+e^{-c_{1}}C_{1}(A_{2})  $$
$$F(\psi) = (1-e^{-c_{1}})F_{1}(A_{1})+e^{-c_{1}}F_{1}(A_{2})  $$
$$H(\psi) = (1-e^{-c_{1}})H_{1}(A_{1})+e^{-c_{1}}H_{1}(A_{2}).  $$
If $$v(a) = c_{1}\delta_{A_{1}}(a)+c_{2}\delta_{A_{2}}(a)+c_{3}\delta_{A_{3}}(a) $$
then in similar manner
$$ \psi(a) = (1-e^{-c_{1}})\delta_{A_{1}}(a)+(e^{-c_{1}}-e^{-c_{1}-c_{2}}) \delta_{A_{2}}(a) +(e^{-c_{1}-c_{2}}-e^{-c_{1}-c_{2}}-c_{2}) \delta_{A_{3}}(a). $$
Three age policy will satisfy the condition $Q(\psi)=1$, so $c_{3} = \infty$. Therefore, we have
$$ \psi(a) = (1-e^{-c_{1}})\delta_{A_{1}}(a)+e^{-c_{1}}(1-e^{-c_{2}})\delta_{A_{2}}(a)+e^{-c_{1}-c_{2}} C_{1}(A_{3}). $$
So, based on the above analysis we can state the following result:

\begin{Proposition}
There exist an optimal vaccination strategy such that individuals in atmost three age classes are vaccinated i.e.
$$ v(a) = c_{1}\delta_{A_{1}}(a)+c_{2}\delta_{A_{2}}(a)+c_{3}\delta_{A_{3}}(a) $$
where we assume that $0 \le A_{1} < A_{2} < A_{3} < \infty$. Also if the constants $c_{1},c_{2}$ and $c_{3}$ are all different form zero then at age $A_{3}$, all the individuals who left at age $A_{1}$ and $A_{2}$ will be vaccinated.
\end{Proposition}

\mysection{Discussion}
It is clear that age plays a crucial role in SARS diseases and other infectious diseases, so introduction of age structure in SIR models  allows us to take care of more realistic scenarios related to disease progression. We propose an age structured SEQIR model in which a proportion of individuals are quarantined (if those individuals are traveling from most affected areas). After the formulation of model, we derive steady state solutions to the model and basic reproduction number is derived by using the appearance of endemic steady state. We construct a Liapunov function to check the stability of disease free equilibrium point. We find the expression for average age of infection and also for constant parameters we derive basic reproduction and average age of infection, we also show that our results matches with SEIR models considered in \cite{MR1814049} if there are no quarantine individuals. With no quarantine individuals we also consider vaccination strategy to our model. We assume that with strict government policies it is possible to recognize those individuals who have already passed through the disease and acquired immunity. So, we vaccinate only susceptible and exposed individuals. We count infected individuals and weigh them with some factor which measure the social impact of one diseased case. With an upper bound on number of infected individuals,  we try to minimize the cost of vaccination. We prove that the vaccination strategy is concentrated on atmost three age classes.   


\end{document}